\documentclass[10pt, a4paper, reqno]{amsart}
\usepackage{preambolo}
\usepackage{tikz-cd}

\begin{document}


\setcounter{secnumdepth}{3}

\setcounter{tocdepth}{2}

\title{\textbf{Kodaira dimension of almost complex $4$-manifolds with torsion first Chern class}
}

\author[Lorenzo Sillari]{Lorenzo Sillari}

\address{Lorenzo Sillari: Dipartimento di Scienze Matematiche, Fisiche e Informatiche, Unità di Matematica e Informatica, Università degli Studi di Parma, Parco Area delle Scienze 53/A, 43124, Parma, Italy} 
\curraddr{\scshape{Erd\H{o}s Center, HUN-REN R\'{e}nyi Institute of Mathematics, Re\'{a}ltanoda utca 14, H-1053, Budapest, Hungary}}
\email{lorenzo.sillari@unipr.it, sillari.lorenzo@renyi.hu}

\author[Adriano Tomassini]{Adriano Tomassini}

\address{Adriano Tomassini: Dipartimento di Scienze Matematiche, Fisiche e Informatiche, Unità di Matematica e Informatica, Università degli Studi di Parma, Parco Area delle Scienze 53/A, 43124, Parma, Italy}
\email{adriano.tomassini@unipr.it}

\maketitle

\begin{abstract} 
\noindent \textsc{Abstract}. In this paper we investigate the Kodaira dimension of almost complex $4$-manifolds with torsion first Chern class. First, we prove that, if the almost complex structure is also tamed, the only possible values for the Kodaira dimension are $0$ or $-\infty$. This is done by developing the theory of pseudoholomorphic structures on vector bundles. In arbitrary dimension, we study infinitesimal deformations of structures with pseudoholomorphically torsion canonical bundle. We compute their tangent space and, under suitable assumptions we prove an unobstructedness theorem in the spirit of Bogomolov--Tian--Todorov. Together, our results allow to fully describe non-integrable infinitesimal deformations of complex structures on $K3$ and Enriques surfaces in terms of their Kodaira dimension.

\end{abstract}

\blfootnote{  \hspace{-0.55cm} 
{\scriptsize 2020 \textit{Mathematics Subject Classification}. Primary: 32Q60, 32L05; Secondary: 53C15, 32G07. \\ 
\textit{Keywords: almost complex Kodaira dimension, almost complex manifold, canonical bundle, deformations of almost complex structures, pseudoholomorphic section, pseudoholomorphic vector bundle.}
\vspace{.1cm}

\noindent The first author is partially supported by the GNSAGA project ``Progetti di Ricerca 2025 - CUP E53C24001950001". The authors are partially supported by GNSAGA of INdAM and by the Project PRIN 2022 “Real and Complex Manifolds: Geometry and Holomorphic Dynamics” (code 2022AP8HZ9)}}

\section{Introduction}\label{sec:intro}

Enriques--Kodaira classification of complex surfaces is a cornerstone of complex geometry. For a complex surface $(M,J)$, the classification establishes a deep connection between the underlying topology of $M$ and key invariants of the complex structure $J$, such as the Kodaira dimension, the plurigenera, and the first Chern class. Among those, the Kodaira dimension $\kappa_J$, taking values in $\{-\infty, 0,1,2\}$, measures the growth rate of the number of holomorphic sections of tensor powers of the canonical bundle, allowing one to separate complex surfaces into four classes.

Although $\kappa_J$ is defined for complex manifolds of any dimension, its properties are most distinctive in complex dimension two. Here, it is not only a smooth invariant \cite{FQ94} but also determines the geometric structure of $(M,J)$. In addition, imposing constraints on $J$ restricts the possible values of $\kappa_J$. A fundamental one is given by the first Chern class $c_1(J) \in H^2(M;\Z)$, which coincides (up to a sign) with the Chern class of the canonical bundle. A particularly important case is when $c_1 = 0$, which corresponds to the canonical bundle being trivial as a smooth vector bundle. In this case, $\kappa_J$ is either $0$ or $-\infty$, see, e.g., \cite{Bau08}, depending on whether or not (a tensor power of) the canonical bundle is holomorphically trivial.

The efficacy of $\kappa_J$ in the classification inspires the search for analogous invariants for other geometric structures. Indeed, a definition of Kodaira dimension for symplectic $4$-manifolds was given in \cite{Li06a}, see also \cite{Mcd90}, and has proven to be particularly fruitful in symplectic geometry \cite{Li06}. In higher dimensions, finding an appropriate definition is still an open problem, cf.\ \cite{Li19}.

A significant recent development is the extension of the Kodaira dimension to almost complex manifolds of arbitrary dimension \cite{CZ23}. Mirroring the pattern seen in the complex and symplectic cases, it exhibits a particular behavior in dimension four, where it can be used to establish relations with the underlying geometry of $(M,J)$, see \cite{CZ23} and \cite{CZ24}.
\vspace{.2cm}

This paper studies the Kodaira dimension of almost complex $4$-manifolds with vanishing or torsion first Chern class. To achieve our goal, we need to focus first on pseudoholomorphic structures on complex line bundles over almost complex manifolds, building on the theory developed in \cite{CZ23}, \cite{BT96}, \cite{Kru07} and \cite{LS01}. The key notion is the Iitaka dimension of pseudoholomorphic structures on line bundles \cite{CZ23}, which generalizes the classical invariant for holomorphic line bundles introduced in \cite{Iit71}. Our first main result is the following.

\begin{mythm}{\ref{thm:iitaka}}
    Let $(M,J)$ be a compact tamed almost complex $4$-manifold. Let $L \longrightarrow (M,J)$ be a complex line bundle with torsion $c_1 (L)$, and let $\bar \partial_L$ be a pseudoholomorphic structure on $L$. Then the Iitaka dimension of $(L, \bar \partial_L)$ is either $0$ or $-\infty$.
\end{mythm}

As an immediate consequence, we show that the Kodaira dimension of an almost complex $4$-manifold with torsion $c_1$ satisfying the tameness assumption can only be $0$ or $-\infty$, and we precisely characterize when it vanishes.

\begin{mycor}{\ref{cor:kodaira:4}}
    The Kodaira dimension of a compact tamed almost complex $4$-manifold with torsion first Chern class is either $0$ or $-\infty$. The Kodaira dimension is $0$ if and only if the canonical bundle is pseudoholomorphically torsion.
\end{mycor}

This can be thought as a generalization to the almost complex setting of the already recalled fact that a complex surface with $c_1=0$ has non-positive Kodaira dimension. However, without integrability, the proof becomes highly non-trivial and requires additional assumptions that guarantee non-existence of certain pseudoholomorphic curves, cf.\ Corollary \ref{cor:no:PH}.

Our approach provides a new method to establish constraints on the Kodaira dimension of almost complex $4$-manifolds without resorting to explicit computations. This complements Theorem 6.7 in \cite{CZ23}, which deals with the case $\kappa_J=1$. Previously, determining $\kappa_J$ typically required either to solve PDEs or to show that they have no solution, see \cite{CNT24}, \cite{CNT20}, and \cite{CNT21}.
\vspace{.2cm}

By Corollary \ref{cor:kodaira:4}, we know that on certain manifolds the Kodaira dimension is either $0$ or $-\infty$, and that both values can coexist on the same connected component of the space of almost complex structures, see, e.g., \cite{ST25}. Since the case $\kappa_J=-\infty$ corresponds to the generic situation (cf.\ Section 3 in \cite{Bry06}), the next step is to understand the space of almost complex structures with vanishing Kodaira dimension.

We can attack the problem by studying the space of pseudoholomorphically trivial line bundles $\mathcal{D}_{PH}$ over a fixed almost complex manifold (of arbitrary dimension). We consider the action of the identity component of the diffeomorphism group on $\mathcal{D}_{PH}$ to get a moduli space $\mathcal{M}$, and we compute its tangent space at $J_0$, where $J_0$ satisfies extra assumptions, including integrability. Namely, if $\sigma_0$ is a pseudoholomorphic trivialization of the canonical bundle, we give a description of the tangent space to $\mathcal{M}$ at $(J_0, \sigma_0)$ in terms of $\mathcal{L}$, the tangent space to the space of all almost complex structures, and the action of the operators $\eth$ and $\bar \partial$, suitable extensions to $\mathcal{L}$ of the operators $\partial$ and $\bar \partial$. An important role in obtaining a nice description in terms of differential operators is played by a mild form of $\partial \bar \partial$-lemma.

\begin{mythm}{\ref{thm:tangent:space}}
    Let $M$ be a compact almost complex $2m$-manifold. Let $(J_0, \sigma_0) \in \mathcal{D}_{PH}$ and consider the space
    \[
    \mathcal{M} \coloneqq  \mathcal{D}_{PH} / \Diff_0(M).
    \]
    If $J_0$ is integrable and it satisfies the $\partial \bar \partial$-lemma on $(m,1)$-forms, then we have that
    \[
    T_{(J_0, \sigma_0)} \mathcal{M} = \frac{\mathcal{L} \cap \ker \eth   \bar \partial}{\mathcal{L} \cap \Ima \bar \partial}.
    \]  
\end{mythm}

The result naturally extends to the space $\mathcal{M}^k_{tor}$ of almost complex structures with pseudoholomorphically torsion canonical bundle, see Theorem \ref{thm:PH:torsion}. To complete our analysis, we study the local deformation theory of these moduli spaces, and we establish that, under appropriate conditions, the moduli space of almost complex structures with pseudoholomorphically torsion canonical bundle is unobstructed.

\begin{mythm}{\ref{thm:unobstructed}}
    Let $M$ be a compact almost complex $2m$-manifold. Let $J_0$ be a fixed almost complex structure on $M$ with pseudoholomorphically torsion canonical bundle. Suppose that $J_0$ is integrable, that $H^{0,1}_{\bar \partial } =\{0\}$, and that the $\partial  \bar \partial $-lemma holds. Then the moduli space $\mathcal{M}^k_{tor}$ is unobstructed in a neighborhood of $J_0$.
\end{mythm}

Theorem \ref{thm:unobstructed} should be thought as a broad generalization of  Bogomolov--Tian--Todorov's Theorem \cite{Bog78, Tia87, Tod89}. Indeed, if $J_0$ is a Calabi--Yau structure and we restrict to the moduli space of integrable deformations, our result recovers the classical one.

The first Chern class being torsion plays, again, a fundamental role, as it provides the link between results about the canonical bundle and conclusions about the Kodaira dimension. By combining Corollary \ref{cor:kodaira:4} and Theorem \ref{thm:unobstructed}, we get a complete local description of the moduli space of almost complex structures with vanishing Kodaira dimension on K\"ahler complex surfaces with torsion canonical bundle, namely on $K3$ and Enriques surfaces.

\begin{mythm}{\ref{thm:K3}}
    Let $(M,J_0)$ be a $K3$ or Enriques surface and let $\mathcal{K}_0(M)$ be the space of almost complex structures on $M$ with vanishing Kodaira dimension. Then the moduli space $\mathcal{K}_0(M)/\Diff_0(M)$ is unobstructed in a neighborhood of $J_0$. Moreover, small deformations of $J_0$ with vanishing Kodaira dimension are in bijection with classes
    \[
    [L] \in  \frac{\mathcal{L} \cap \ker \eth   \bar \partial}{\mathcal{L} \cap \Ima \bar \partial}.
    \]
\end{mythm}

To the best of our knowledge, this is the first known characterization of the space of almost complex structures with Kodaira dimension $0$. This should be compared to the fact that integrable deformations of $K3$ and Enriques surfaces all have vanishing Kodaira dimension.
\vspace{.2cm}

The paper is structured as follows. Section \ref{sec:phvb} deals with the basics on pseudoholomorphic vector bundles and line bundles. In Section \ref{sec:iitaka}, we prove our main result on the Iitaka dimension and the Kodaira dimension of tamed almost complex $4$-manifolds. Finally, Section \ref{sec:deformations} deals with infinitesimal deformations of structures with pseudoholomorphically torsion canonical bundle, and the consequences on $K3$ and Enriques surfaces.

\subsection*{Acknowledgments.}
The first author wishes to thank the Erd\H{o}s Center and the R\'enyi Institute, where part of the present paper was written, for the pleasant and stimulating work environment. The authors are grateful to Tian-Jun Li and Riccardo Piovani for useful comments on an early version of this paper.

\section{Pseudoholomorphic vector bundles}\label{sec:phvb}

The goal of this section is to collect basic properties of pseudoholomorphic vector bundles over almost complex manifolds, and to introduce the Iitaka dimension of pseudoholomorphic line bundles. Aside from some original result, our main references for the general theory of pseudoholomorphic vector bundles are \cite{BT96}, \cite{Kru07} and \cite{LS01}. For the approach leading to the definition of Kodaira dimension of almost complex manifolds, see \cite{CZ23} and \cite{CZ24}.

\subsection{Preliminaries and notation.}\label{sec:notation}

For brevity, when referring to vector bundles, sections, or functions, we will often write \textit{PH} instead of \textit{pseudoholomorphic}. 

Let $(M, J)$ and $(N,J')$ be two almost complex manifolds. A map $f \colon (M, J) \rightarrow (N,J')$ is said to be \emph{$(J, J')$-holomorphic} if 
\[
df \circ J = J' \circ df.
\]
When there is no ambiguity on the almost complex structures on the source and on the target manifolds, we will say, with a slight abuse of terminology, that $f$ is \emph{pseudoholomorphic}. A submanifold $\Sigma$ of $(M,J)$ is a \textit{PH submanifold} if $T \Sigma$ is preserved under the action of $J$. 

We denote by $\Lambda^k(M)$ the bundle of $k$-forms on $M$ and by
\[
\Lambda^k(M) = \bigoplus_{p+q=k} \Lambda^{p,q}(M)
\]
the bigraded decomposition into $(p,q)$-forms induced on it by $J$. The corresponding spaces of smooth sections are $A^k(M)$ and $A^{p,q}(M)$. We denote by $\Diff_0(M)$ the group of diffeomorphisms of $M$ that are isotopic to the identity.

Recall that an almost complex manifold $(M,J)$ is said to be \emph{tamed} if there exists a symplectic form $\omega \in A^2(M)$ taming $J$, that is, such that $\omega_x (v, Jv) >0$ for every tangent vector $v \in T_xM$, and that $J$ is said to be \emph{integrable} if its Nijenhuis tensor $N_J$ vanishes, which induces on $M$ the structure of a complex manifold. A complex manifold is said to satisfy the $\partial \bar \partial$-lemma on $(p,q)$-forms if 
\[
A^{p,q}(M) \cap \ker \partial \cap \Ima \bar \partial = A^{p,q}(M) \cap \ker \bar \partial \cap \Ima \partial = A^{p,q}(M) \cap \Ima \partial \bar \partial.
\]
If the $\partial \bar \partial$-lemma holds for every $(p,q)$, we simply say that $(M,J)$ satisfies the $\partial \bar \partial$-lemma.

In the complex case, a holomorphic vector bundle is a real vector bundle $\pi \colon E \rightarrow M$ where both $E$ and $M$ are complex manifolds and the projection $\pi$ is holomorphic. This naturally induces on $E$ the structure of a complex vector bundle.

\subsection{Pseudoholomorphic vector bundles.}

The following generalization of the notion of holomorphic vector bundle to the almost complex setting is due to Lempert and Sz\H{o}ke, see Definition 4.1 in \cite{LS01}.

\begin{definition}\label{def:phvb}
    A \emph{pseudoholomorphic vector bundle} is a real vector bundle
    \[
    \pi \colon (E, \J) \longrightarrow (M,J),
    \]
    where $(E,\J)$ and $(M,J)$ are almost complex manifolds, such that:
    \begin{itemize}
        \item[(i)] the projection $\pi$ is $(\J, J)$-holomorphic;
        \item[(ii)] for all $x \in M$, the fiberwise sum $E_x \times E_x \rightarrow E_x$ is $(\J \times \J, \J)$-holomorphic.
    \end{itemize}
\end{definition}

An alternative definition is given by Kruglikov, see Definition 1 in \cite{Kru07}, by replacing condition (ii) with the following differential-geometric condition:
    \begin{itemize}
        \item[(ii')] for all $x \in M$, the structure $\J$ has constant coefficients on $E_x$ and there exists a connection $\nabla$ on $E$ providing a split short exact sequence
        \[
        0 \rightarrow E_x \rightarrow T_yE \rightarrow T_xM \rightarrow 0, 
        \]
        for $y \in E$ and $x = \pi(y)$.
    \end{itemize}

By Proposition 19 in \cite{Kru07}, the two definitions are equivalent. In addition, a PH vector bundle naturally inherits the structure of a complex vector bundle induced by fiberwise sum, hence a well-defined fiberwise multiplication $\C \times E_x \rightarrow E_x$, that in general is not PH, cf.\ \cite{Kru07} and \cite{LS01}. 

Let $N'_\J$ be the Nijenhuis tensor of $\J$ computed on the horizontal bundle (with respect to the connection $\nabla$ associated to $\J$) in the first argument and on the vertical bundle in the second argument. By Proposition 13 in \cite{Kru07}, fiberwise multiplication is PH if and only if $N_\J' =0$. PH vector bundles with $N_\J'=0$ are called \textit{normally integrable}, and they coincide with \textit{bundle almost complex structures}, as defined by de Bartolomeis and Tian, see Definition 1.2 in \cite{BT96}.

\begin{definition}\label{def:bacs}
     A \emph{bundle almost complex structure} is a real vector bundle
    \[
    \pi \colon (E, \J) \longrightarrow (M,J),
    \]
    where $(E,\J)$ and $(M,J)$ are almost complex manifolds, such that:
    \begin{itemize}
        \item[(i)] the projection $\pi$ is $(\J, J)$-holomorphic;
        \item[(ii)] for all $x \in M$, the fiberwise sum $E_x \times E_x \rightarrow E_x$ is $(\J \times \J, \J)$-holomorphic;
        \item [(iii)] for all $x \in M$, the fiberwise multiplication $\C \times E_x \rightarrow E_x$ is $(J_2 \times \J, \J)$-holomorphic, where $J_2$ is the standard complex structure on $\C$.
    \end{itemize}
\end{definition}

\begin{remark}\label{remark:multiplication}

If $\J$ and $J$ are both integrable, one recovers the usual notion of holomorphic vector bundle. Under the integrability assumption, fiberwise sum and multiplication by a complex number are always holomorphic maps, and there is no need to ask for condition \emph{(ii)} for PH vector bundles or \emph{(ii)} and \emph{(iii)} for normally integrable PH vector bundles. 
\end{remark}

\begin{remark}
Let $(E, \J) \rightarrow (M,J)$ be a normally integrable PH vector bundle of rank $r$, and fix $x \in M$. Consider the complex structure on the fibers $E_x \cong \C^r$ induced by the PH vector bundle structure. Then $\J$ restricted to the fiber $E_x$ coincides with the standard complex structure on $\C^r$, cf.\ Definition 1.2 in \cite{BT96}. 

\end{remark}

The choice of a connection is essential when developing differential geometry of PH vector bundles, since it allows us to work PH structures, which generalize Dolbeault operators.

\begin{definition}
    Let $E \rightarrow (M,J)$ be a complex vector bundle, where $(M,J)$ is an almost complex manifold. A \emph{pseudoholomorphic structure} on $E$ is an operator
    \[
    \bar \partial_E \colon \Gamma (E) \longrightarrow \Gamma (E \otimes \Lambda^{0,1} (M))
    \]
    satisfying the Leibniz rule
    \[
    \bar \partial_E (f \sigma) = \bar \partial f \otimes \sigma + f \bar \partial_E \sigma,
    \]
    for all $f \in C^\infty (M)$ and $\sigma \in \Gamma (E)$, where $\bar \partial$ is the standard Dolbeault operator on $M$. A section $\sigma \in \Gamma (E)$ is said to be \emph{pseudoholomorphic} if $\bar \partial_E \sigma =0$.
\end{definition}

\begin{remark}

As proved in Section 3 in \cite{CZ23}, PH structures are compatible (eventually after fixing a compatible Hermitian metric on $E$) with the standard operations on vector bundles of duality, direct sum, and tensor product.

\end{remark}

\begin{remark}\label{rem:connection}
By Proposition 1.3 in \cite{BT96}, every PH structure $(E, \bar \partial_E)$ induces on $E$ a structure of normally integrable PH vector bundle, given by the formula
\begin{equation}\label{eq:BT:onetoone}
    df (\J X) = -i (2 \bar \partial_E f - df) (X),
\end{equation}
where $X$ is a vector field on $E$ and $f \in C^\infty (E)$. Conversely, if $(E,\J)$ is a PH vector bundle, then $E$ naturally admits a PH structure $\bar \partial_E$. For more details on the relation between the operator $\bar \partial_E$ and the connection $\nabla$ given in condition (ii'), we refer the reader to Section 1 in \cite{BT96}.
    
\end{remark}

In what follows, we will heavily exploit the equivalence of the differential geometric approach, based on the operator $\bar \partial_E$, and the topological approach, based on PH submanifolds with respect to the almost complex structures $\J$ and $J$. In particular, we recall that a section $\sigma \in \Gamma(E)$ is PH with respect to a PH structure $\bar \partial_E$ if and only if $d \sigma \circ J = \J \circ d \sigma$, where $\J$ is the almost complex structure induced on $E$ by $\bar \partial_E$ as in \eqref{eq:BT:onetoone}. For more details on the topic, we refer to Lemma 4.1 in \cite{CZ23} or to \cite{CZ24}.

\begin{definition}\label{def:equivalence}
Let $(E, \J_E)$ and $(F, \J_F)$ be two PH vector bundles in the sense of Definition \ref{def:phvb}. We say that they are \textit{PH equivalent} if there exists an isomorphism $\varphi \colon E \rightarrow F$ of real vector bundles such that $d \varphi \circ \J_E = \J_F \circ d \varphi$. 

Let $(E, \bar \partial_E)$ and $(F, \bar \partial_F)$ be two complex vector bundles $E$ and $F$ over $(M,J)$ endowed with two PH structures $\bar \partial_E$ and $\bar \partial_F$. We say that they are \textit{PH equivalent} if there exists an isomorphism $\varphi \colon E \rightarrow F$ of real vector bundles such that the following diagram commutes
\[
\begin{tikzcd}
    \Gamma(E) \arrow[dd, "\bar \partial_E"'] \arrow[r, "\varphi"]          & \Gamma(F) \arrow[dd, "\bar \partial_F"]   \\
     & \\
    \Gamma(E) \otimes A^{0,1}(M) \arrow[r, "\varphi \otimes (\varphi|_M^{-1})^*"'] & \Gamma(F) \otimes A^{0,1}(M)
\end{tikzcd}
\]
\end{definition}

Observe that if $(E, \J_E)$ and $(F, \J_F)$ are PH equivalent, then one of them is normally integrable if and only if the other one is.

\begin{lemma}
    Two PH structures $(E, \bar \partial_E)$ and $(F, \bar \partial_F)$ are PH equivalent if and only if the corresponding normally integrable PH vector bundles $(E, \J_E)$ and $(F, \J_F)$ are PH equivalent.
\end{lemma}
\begin{proof}
    Let $(E, \J_E)$ and $(F, \J_F)$ be two PH vector bundles over $(M,J)$. Let $\varphi \colon E \rightarrow F$ be a vector bundle isomorphism. The corresponding PH structures $(E, \bar \partial_E)$ and $(F, \bar \partial_F)$ are PH equivalent if and only if for every $\sigma \in \Gamma (E)$ we have that
    \begin{equation}\label{eq:PHstructure:iso}
    \bar \partial_F (\varphi \circ \sigma) = \varphi \otimes (\varphi|_M^{-1})^* ( \bar \partial_E \sigma).
    \end{equation}
    Fix a local trivialization $\{ s_j \}_{j=1}^r$ of $E$ such that $\sigma = f_j s_j$, where $r$ is the complex rank of $E$, $f_j \in C^\infty(M)$, and we are summing over repeated indices. Let $^E\theta_j^k$ be the connection $1$-form of $\bar \partial_E$ and $^F\theta_j^k$ those of $\bar \partial_F$, cf.\ Remark \ref{rem:connection} or \cite{BT96}. We immediately see that
    \[
    \bar \partial_F (\varphi \circ s_j) = (\varphi \circ s_k) \otimes  \, ^F\theta_j^k.
    \]
    On the other side, we also have that
    \[
    \varphi \otimes (\varphi|_M^{-1})^* (\bar \partial_E s_j )= \varphi \otimes (\varphi|_M^{-1})^* ( s_k \otimes\,  ^E\theta_j^k) = (\varphi \circ s_k) \otimes (\varphi|_M^{-1})^* (\, ^E\theta_j^k)
    \]
    Since $\bar \partial_E$ and $\bar \partial_F$ are $\C$-linear and act as derivations, we can rewrite \eqref{eq:PHstructure:iso} for $\sigma = s_j$ as
    \[
    (\varphi \circ s_k) \otimes  \, ^F\theta_j^k = (\varphi \circ s_k) \otimes (\varphi|_M^{-1})^* (\, ^E\theta_j^k),
    \]
    in order to get the condition
    \begin{equation}\label{eq:connection:forms}
        ^E\theta_j^k = (\varphi|_M)^* (\, ^F\theta_j^k),
    \end{equation}
    which is equivalent to \eqref{eq:PHstructure:iso}. We now use the explicit expression of $\J_E$ and $\J_F$ in terms of $\bar \partial_E$ and $\bar \partial_F$ given by \eqref{eq:BT:onetoone}, see also the first equation on page 74 in \cite{LS01}. For every smooth function $f \colon F \rightarrow \C$, $g = f \circ \varphi \colon E \rightarrow \C$, and every vector field $X \in \Gamma (TE)$, we have that
    \begin{equation}\label{eq:equiv:1}
    df \circ d\varphi \circ \J_E (X) = dg (\J_E X) = -i (2 \bar \partial_E g -dg)(X).
    \end{equation}
    Similarly, by the expression of $\J_F$ and by writing an arbitrary vector field on $F$ as $d \varphi (X)$, we get
    \begin{equation}\label{eq:equiv:2}
    df \circ \J_F \circ d\varphi (X) = - i (2 \bar \partial_F f - df) \circ d \varphi (X).
    \end{equation}
    If \eqref{eq:PHstructure:iso} holds, equation \eqref{eq:connection:forms} gives that $\bar \partial_F f \circ d\varphi = \bar \partial_E (f \circ \varphi)$, hence the two right hand sides of \eqref{eq:equiv:1} and \eqref{eq:equiv:2} coincide. This readily implies that $d \varphi \circ \J_E = \J_F \circ d \varphi$. The converse implication follows by comparing \eqref{eq:equiv:1} and \eqref{eq:equiv:2}, and by the arbitrariety of $f$ and $X$.
\end{proof}

\begin{definition}\label{def:triviality}
A PH vector bundle $(E, \J_E)$ of complex rank $r$ is said to be \textit{PH trivial} if it is PH equivalent to the bundle $M \times \C^r$ endowed with the almost complex structure acting as $J$ on $M$ and as the standard complex structure on $\C^r$.
\end{definition}

\begin{lemma}\label{lemma:trivial:equivalence}
    A normally integrable PH vector bundle of complex rank $r$ is PH trivial if and only if it admits $r$ independent PH sections.
\end{lemma}
\begin{proof}
    Let $J_{2r}$ be the standard complex structure on $\C^r$. The trivial bundle $M \times \C^r$ endowed with the almost complex structure $\J \coloneqq J \times J_{2r}$ is normally integrable and it admits $r$ independent sections
    \begin{align*}
    s_j \colon M &\longrightarrow M \times \C^r\\
    x &\longmapsto (x, e_j),
    \end{align*}
    where $e_j$, for $j =1, \ldots, r$, is the standard unitary basis of $\C^r$. To see that each $s_j$ is PH, one can either verify that $\bar \partial s_j=0$, since $s_j$ has constant coefficients, or to directly check that 
    \[
    d s_j \circ J = \J \circ ds_j.
    \]
    Similarly, if $(E, \J_E) \rightarrow (M,J)$ is PH equivalent to $M \times \C^r$ via the map $\varphi \colon E \rightarrow M \times \C^r$, then we can find $r$ independent PH sections of $E$ by pulling back the $s_j$ via $\varphi$.
    
    Conversely, let $(E, \J_E) \rightarrow (M,J)$ be a normally integrable PH vector bundle of rank $r$ admitting $r$ independent PH sections. The map
    \begin{align*}
        \varphi \colon M \times \C^r &\longrightarrow E,\\
        (x, v = \sum_{j=1}^r v_j e_j) &\longmapsto (x, \sum_{j=1}^r v_j s_j(x))
    \end{align*}
    provides a smooth trivialization of $E$. Since fiberwise sum and multiplication are PH maps, the map $\varphi$ is PH, proving the claim. Alternatively, to directly check that the trivialization is PH, we show that $d \varphi \circ \J = \J_E \circ d \varphi$. Explicitly, let $(X,Y) \in T_{(x,v)} (M \times \C^r)$. The action of $d \varphi$ is
    \[
    d \varphi_{(x,v)} (X,Y) = \sum_{j=1}^r(v_j ds_j (X) + Y_j s_j(x)),
    \]
    where $Y = \sum\limits_j Y_j e_j \in T_v \C^r$. By applying $\J_E$ to both sides, we get the equality 
    \begin{equation}\label{eq:trivial:1}
    \J_E \circ d \varphi_{(x,v)} (X,Y) = \sum_{j=1}^r(v_j \J_E \circ ds_j (X) + i Y_j s_j(x)),
    \end{equation}
    where we used that $\J_E$ acts as multiplication by $i$ on the fibers of $E$. On the other side, we also have that 
    \begin{equation}\label{eq:trivial:2} 
    d \varphi_{(x,v)} \circ \J (X,Y) = d \varphi_{(x,v)} (JX,iY) = \sum_{j=1}^r(v_j ds_j (JX) + i Y_j s_j(x)).
    \end{equation}
    Since the sections $s_j$ are PH, the relation $ds_j \circ J = \J_E \circ ds_j$ holds, showing that \eqref{eq:trivial:1} and \eqref{eq:trivial:2} coincide, and concluding the proof.
\end{proof}

\subsection{Line bundles and their Iitaka dimension.}
Let $L$ be a complex line bundle over $M$ and let  $\bar \partial_L$ be a PH structure on it. Tensor powers of $L$ inherit a PH structure that we still denote by $\bar \partial_L$. Explicitly, for $l \ge 2$, we have that $(L^{\otimes l}, \bar \partial_L)$ is a PH line bundle with $\bar \partial_L$ defined recursively by the formula
\[
\bar \partial_L (\sigma_1 \otimes \sigma_2 ) \coloneqq \bar \partial_L \sigma_1 \otimes  \sigma_2 + \sigma_1 \otimes \bar \partial_L \sigma_2,
\]
where $l = l_1 +l_2$, $\sigma_1 \in \Gamma (L^{\otimes l_1})$, and $\sigma_2 \in \Gamma (L^{\otimes l_2})$. Denote by $P_l$ the (complex) dimension of the space of PH sections of $L^{\otimes l}$. The number $P_l$ depends on the choice of line bundle $L$ and on the choice of PH structure $\bar \partial_L$.

The \textit{Iitaka dimension of $(L, \bar \partial_L)$}, see Definition 4.2 in \cite{CZ23}, is the number 
\[
\kappa (M,L, \bar \partial_L) \coloneqq
\begin{cases}
    - \infty & \text{if $P_l =0$ for all $l \ge 1$,}\\
    \limsup\limits_l \frac{\log P_l}{\log l} & \text{otherwise.}
\end{cases}
\]

Let $K_J \coloneqq \Lambda^{m,0}(M)$ be the bundle of $(m,0)$-forms on $(M,J)$. We call $K_J$ the \textit{canonical bundle of $(M,J)$}. The operator
\[
\bar \partial \colon \Lambda^{m,0}(M) \longrightarrow \Lambda^{m,1}(M)
\]
induces a PH structure on $K_J$, hence on its tensor powers $K_J^{\otimes l}$ for $l \ge 1$, that we still denote by $\bar \partial$. The \textit{Kodaira dimension of $(M,J)$} is the Iitaka dimension of $(K_J, \bar \partial)$. Explicitly, we have that 
\[
\kappa_J \coloneqq
\begin{cases}
    - \infty & \text{if $P_l =0$ for all $l \ge 1$,}\\
    \limsup\limits_l \frac{\log P_l}{\log l} & \text{otherwise.}
\end{cases}
\]
The choice of notation $\kappa_J$ underlines that the Kodaira dimension depends only on the choice of almost complex structure on $M$. Moreover, if $M$ is compact then the $P_l$ are finite numbers and $\kappa_J < +\infty$ \cite{CZ23}. If $J$ is integrable, the definition of $\kappa_J$ coincides with the classical definition of Kodaira dimension for complex manifolds.

\section{Iitaka dimension of PH structures of line bundles with \texorpdfstring{$c_1 =0$}{}}\label{sec:iitaka}

In this section we compute the Iitaka dimension of PH structures of complex line bundles with torsion $c_1$ over tamed almost complex $4$-manifolds. As an application, we are able to determine the Kodaira dimension for a wide class of almost complex structures.
\vspace{.2cm}

Let $(M,J)$ be a compact almost complex $4$-manifold and let $L \longrightarrow (M,J)$ be a complex line bundle. Suppose that $c_1(L) \in H^2 (M;\Z)$ is torsion, so that a certain tensor power of $L$ is smoothly trivial. If we fix a PH structure $\bar \partial_L$, the Iitaka dimension of $(L, \bar \partial_L)$ is related to whether or not $L$ is PH trivial, or PH torsion (any of its tensor powers is PH trivial). This is made precise in the following theorem, where we see that the condition of $c_1$ being torsion imposes severe restrictions on the Iitaka dimension of PH structures on $L$.

\begin{theorem}\label{thm:iitaka}
    Let $(M,J)$ be a compact tamed almost complex $4$-manifold. Let $L \longrightarrow (M,J)$ be a complex line bundle with torsion $c_1 (L)$, and let $\bar \partial_L$ be a PH structure on $L$. Then we have that
    \[
    \kappa (M, L, \bar \partial_L) \in \{ -\infty, 0 \}.
    \]
\end{theorem}
\begin{proof}
    Let $\bar \partial_L$ be an arbitrary PH structure on $L$ and let $\J$ be the corresponding normally integrable almost complex structure on the total space of $L$. Observe that, by condition (i) in Definition \ref{def:phvb}, the structure $\J$ restricted to $TM$ coincides with $J$. 
    
    If $L^{\otimes k}$ has no PH sections for each $k$, then its Iitaka dimension is $-\infty$ by definition. 
    Suppose, instead, that there exists at least one PH section $\hat{\sigma}$ of $L^{\otimes k}$ (endowed with the natural PH structure $\bar \partial_L$) for some $k$, and let $\tau$ be the torsion order of $c_1(L)$. Then, the bundle $L^{\otimes \tau k}$ admits a PH section $\sigma = \hat{\sigma}^{\otimes \tau}$. We claim that $\sigma$ must be never-vanishing. Indeed, by Proposition 4.1 in \cite{Zha18}, if $\sigma$ is a PH section of a PH line bundle with $\J|_{TM} = J$, its zero set supports a PH $1$-subvariety $\Sigma$ of $(M,J)$ whose homology class is $PD (c_1 (L^{\otimes \tau k})) = PD (\tau k c_1 (L)) = [0] \in H_2(M;\Z)$. By the tameness assumption, such a homology class cannot exist (since there are no non-constant homologically trivial PH curves on tamed almost complex manifolds), hence $\sigma$ is never-vanishing.

    To complete the proof, we show that if $L^{\otimes k}$ admits a never-vanishing PH section, then the Iitaka dimension of $(L, \bar \partial_L)$ must be $0$. The proof is a generalization of that of Lemma 2.2 in \cite{ST25a}. Let $k_0$ be the minimum integer such that $L^{\otimes k_0}$ admits a never-vanishing PH section, and let $\sigma$ be such a section. We show that the dimension of the space of PH sections of $L^{\otimes k}$ is $1$ if $k$ is a multiple of $k_0$ and $0$ otherwise. This immediately implies that $\kappa (M,L,\bar \partial_L) =0$.
    
    Let $g \sigma$ be another PH section of $L^{\otimes k_0}$, where $g \in C^\infty (M)$ is not necessarily never-vanishing. Then we have that
    \[
    0 = \bar \partial_L (g  \sigma) = \bar \partial g \otimes  \sigma,
    \]
    where $\bar \partial$ is the Dolbeault operator on $M$. Since $ \sigma$ is never-vanishing, this implies that $\bar \partial g =0$ and, by compactness of $M$, that $g$ is constant. Hence, the space of PH sections of $L^{\otimes k_0}$ is one-dimensional. Similarly, if $k = p k_0$, we have that $\sigma^{\otimes p}$ is a never-vanishing PH section of $L^{\otimes k}$ and every other PH section of $L^{\otimes k}$ is a constant multiple of it.
    
    Now suppose, by contradiction, that for $p < k_0$ there exists a (not necessarily never-vanishing) PH section $\eta$. Then $\eta^{\otimes k_0}$ and $ \sigma^{\otimes p}$ are both PH sections of $L^{\otimes p k_0}$. By the same argument of the first part of the proof, we have that $\eta^{\otimes k_0}$ is a constant multiple of $\sigma^{\otimes p}$, which is never-vanishing. In particular, one has that $\eta^{\otimes k_0}$ is never-vanishing, and so is $\eta$. This contradicts the minimality assumption on $k_0$ and proves that $L^{\otimes p}$ has no PH sections for $p < k_0$.
    
    To conclude, let $p$ be any integer that is not a multiple of $k_0$, and let $q$ be the integer such that $q k_0 < p < (q+1) k_0$. Suppose that $\eta$ is a PH section of $L^{\otimes p}$. Then we can write 
    \[
    \eta = \eta' \otimes \sigma^{\otimes q},
    \]
    with $\eta' \in L^{\otimes (p-q k_0)}$, and we get
    \begin{align*}
    0 &= \bar \partial_L ( \eta ) = \bar \partial_L (\eta') \otimes \sigma^{\otimes q},
    \end{align*}
    that implies $\bar \partial_L \eta' =0$, since $\sigma$ is never-vanishing. This would give a PH section of $L^{\otimes (p-q k_0)}$, contradicting again the minimality of $k_0$.
\end{proof}

As a direct consequence, we can compute the Kodaira dimension of tamed almost complex $4$-manifolds with torsion $c_1$, and characterize when they have Kodaira dimension $0$.

\begin{cor}\label{cor:kodaira:4}
    The Kodaira dimension of a compact tamed almost complex $4$-manifold with torsion first Chern class is either $0$ or $-\infty$. The Kodaira dimension is $0$ if and only if the canonical bundle is PH torsion.
\end{cor}
\begin{proof}
    Take $(L, \bar \partial_L)$ to be the canonical bundle $K_J$ of $(M,J)$ endowed with the standard Dolbeault operator. The canonical bundle has torsion first Chern class since
    \[
    c_1 (K_J) = c_1 (\det (T^*M^{1,0})) = c_1 (T^*M^{1,0}) = -c_1 (J), 
    \]
    while the standard Dolbeault operator induces on the total space of $K_J$ an almost complex structure $\J$ such that $\J|_{TM} = J$. Hence, by Theorem \ref{thm:iitaka}, the Iitaka dimension of $(K_J, \bar \partial)$, that is, the Kodaira dimension of $(M,J)$, is either $0$ or $-\infty$. This proves the first part of the corollary. 

    For the second part, it is easy to see that the Kodaira dimension of an almost complex structure with PH torsion canonical bundle has to be $0$ (even without the tameness assumption and in every dimension), cf.\ Lemma 2.2 in \cite{ST25a}. Conversely, suppose that $J$ is an almost complex structure on $M$ such that $c_1 (J)$ is torsion and $\kappa_J =0$. Then, for all $\varepsilon$ and for $l$ large enough, we have the bound $P_l \le l^\varepsilon$, that is, $P_l$ grows slower than any polynomial in $l$. Suppose by contradiction that $P_l \ge 2$ and let $l_0$ be an integer such that $P_{l_0} \ge 2$. Then $P_{t l_0}$ grows at least polynomially in $t$, contradicting the bound $P_{t l_0} \le t^\varepsilon l_0^\varepsilon$. This shows that $P_l \in \{ 0, 1 \}$ for all $l$. Since $\kappa_J \neq -\infty$, there exists $l_0 \ge 1$ such that $P_{l_0} = 1$. Hence, there is exactly one holomorphic section of $K_J^{\otimes k_0}$, which, with the same argument of the proof of Theorem \ref{thm:iitaka}, must be never-vanishing. This shows that $K_J$ is PH torsion and completes the proof.
\end{proof}

Observe that, with the same proof of Theorem \ref{thm:iitaka}, we have the following result.

\begin{cor}\label{cor:no:PH}
    Let $(M,J)$ be a compact almost complex $4$-manifold without homologically trivial PH curves. Let $L \longrightarrow (M,J)$ be a complex line bundle with torsion $c_1 (L)$, and let $\bar \partial_L$ be a PH structure on $L$. Then we have that
    \[
    \kappa (M, L, \bar \partial_L) \in \{ -\infty, 0 \}.
    \]
\end{cor}
\begin{proof}
    The proof follows step by step that of Theorem \ref{thm:iitaka}, with the exception that, without the tameness assumption, we cannot conclude non-existence of homologically trivial PH curves and we have to assume it. 
\end{proof}

\section{Deformations of structures with PH torsion canonical bundle}\label{sec:deformations}

In this section we consider deformations of almost complex structures with PH torsion canonical bundle. We describe the tangent space to such a space of deformations as the kernel of certain differential operators and, under suitable assumptions, we prove an unobstructedness result for the deformation space, in the spirit of Bogomolov--Tian--Todorov's Theorem. As a consequence, we are able to characterize almost complex structures with Kodaira dimension $0$ on $K3$ and Enriques surfaces.
\vspace{.2cm}

To describe small deformations, we adopt a notation similar to that of \cite{BM10} and \cite{BT13}.

Let $M$ be a compact almost complex manifold of real dimension $2m$ and let $\J$ be the space of almost complex structures on $M$. Suppose that $J_0 \in \J$ is an almost complex structure with $c_1 (J_0) =0$. Let $K_0$ be the canonical bundle of $J_0$ and let $\sigma_0$ be a (smooth) trivialization of $K_0$, that is, a never-vanishing $(m,0)$-form on $M$. 

\subsection{Deformation of pairs \texorpdfstring{$(J_0, \sigma_0)$}{}.}
Let $J_t \in \J$, for small $t \in (-\varepsilon, \varepsilon)$, be a smooth curve of almost complex structures passing through $J_0$. Denote by $K_t$ the canonical bundle of $J_t$. Since $c_1$ is invariant under small smooth deformations, we have that $c_1(J_t) =0$ and that $K_t$ is smoothly trivial.

Fix an arbitrary (not necessarily $J_0$-compatible) Riemannian metric $g$ and let $\omega (\cdot, \cdot) \coloneqq -g (\cdot, J_0 \cdot)$ be the corresponding fundamental $2$-form. Then $\omega (\cdot, J_0 \cdot) >0$ and, since strict positivity is an open condition, we still have $\omega (\cdot, J_t \cdot) >0$ for $t$ small enough. By Proposition 1.1.6 in \cite{Aud94}, this implies that $J_t$ can be written as
\[
J_t = (I + L_t) J_0 (I+L_t)^{-1},
\]
where $L_t \in \End (TM)$, with $\lVert L_t \rVert < 1$, $L_t J_0 + J_0 L_t =0$ and $L_t= tL + o(t)$. Note that since $\lVert L_t \rVert < 1$, the endomorphism $(I+L_t)$ is invertible. By taking the time derivative at $t=0$, we see that the tangent space to $\J$ is linearized by
\[
\mathcal{L} = \{ L \in \End (TM): LJ_0 + J_0L =0, \lVert L \rVert < 1\},
\]
which corresponds to the open unit ball in $A^{0,1}(M) \otimes T^{1,0}M$.

The almost complex structure $J_t$ induces a bigraded decompositions of complex vector fields and differential forms that differs from the one induced by $J_0$. More precisely, we have that $X$ is a $(1,0)$-vector field at time $t=0$ if and only if $(I+L_t)X$ is a $(1,0)$-vector field at time $t$. Let $\alpha$ be a smooth $(p,q)$-form and consider the form $\alpha_t = \rho_t \alpha$ defined by
 \begin{multline}\label{eq:rho}
    (\rho_t \alpha) (X_1, \ldots, X_p, Y_1, \ldots, Y_q) :=\\
    \alpha ( (I+L_t)^{-1} X_1, \ldots, (I+L_t)^{-1} X_p,(I+L_t)^{-1} Y_1, \ldots, (I+L_t)^{-1} Y_q),
\end{multline}
where the $X_j$ are $(1,0)$-vector fields and the $Y_j$ are $(0,1)$-vector fields. Then $\alpha$ is a $(p,q)$-form at time $t=0$ if and only if $\alpha_t$ is a $(p,q)$-form at time $t$. As a consequence, the map $\rho_t$ provides an isomorphism 
\[
\rho_t \colon A^{p,q}(M) \xrightarrow{\sim} A^{p,q}_t(M),
\]
where $A^{p,q}_t(M)$ denotes the space of $(p,q)$-forms for $J_t$, that can be extended to an isomorphism between vector valued differential forms as
\begin{align*}
    \rho_t \colon &A^{p,q}(M) \otimes TM^{1,0} \xrightarrow{\sim} A^{p,q}_t(M) \otimes TM^{1,0}_t, \\
    &\alpha \otimes X \longmapsto \rho_t \alpha \otimes (I+L_t)X.
\end{align*}

Using the map $\rho_t$, we can describe trivializations of $K_t$. Indeed, an arbitrary trivialization of $K_t$ is given by 
\[
\sigma_t \coloneqq F_t \rho_t(\sigma_0) \in A^{m,0}_t(M),
\]
where $F_t \in C^\infty(M)$ is never-vanishing. For small $t$ and up to normalization, we can assume that $F_t = 1 + tf + o(t)$, with $f \in C^\infty(M)$. Its linearization is then given by
\[
C^\infty (M)^* \coloneqq \{ f \in C^\infty(M): 1+f \neq 0\}.
\]

Let now $(J_t, \sigma_t)$ be a simultaneous deformation of $J_0$ and of a trivialization of its canonical bundle. By definition, we have that
\[
(J_t, \sigma_t) \in \mathcal{D} = \{ (K, \eta) \in \J \times A^m(M): \eta \in A^{m,0}_K(M) \},
\]
and, by the previous discussion, the tangent space to $\mathcal{D}$ at $(J_0, \sigma_0)$ is
\[
T_{(J_o,\sigma_0)} \mathcal{D} \cong \mathcal{L} \times C^\infty(M)^*.
\]

\subsection{Tangent space to the space of structures with PH trivial canonical bundle.}\label{sec:PH:trivial} Suppose that $K_0$ is PH trivial, that is, that the equation $\bar \partial \sigma_0 =0$ is satisfied, where $\bar \partial$ is the standard $\bar \partial$-operator at time $t=0$. The space of deformations of $(J_0, \sigma_0)$ with PH trivial canonical bundle is
\[
\mathcal{D}_{PH} = \{ (K, \eta) \in \J \times A^m(M): \eta \in A^{m,0}_K(M), \bar \partial_K \eta =0\} \subset \mathcal{D},
\]
where $\bar \partial_K$ is the $\bar \partial$-operator associated to $K$. To describe the tangent space $T_{(J_0, \sigma_0)} \mathcal{D}_{PH}$, we need to introduce several differential operators.

Given $L_t \in \End(TM)$, define a map $\tau_{L_t} \colon A^k (M)\rightarrow A^k(M)$ by setting
\[
\tau_{L_t} (\alpha) (X_1, \ldots, X_k) \coloneqq \sum_{j=1}^k \alpha (X_1, \ldots, L_t X_j, \ldots, X_k).
\]
The \emph{Nijenhuis torsion} \cite{Nij51} of $L_t$ is given, for $X$ and $Y$ vector fields on $M$, by the formula
\[
N_{L_t}(X,Y) = [L_t X,L_t Y]- L_t [L_t X,Y] - L_t [X,L_t Y] +L_t^2 [X,Y].
\]
Then, we can define the operator $\mathcal{N}_{L_t} \colon A^k(M) \rightarrow A^{k+1}(M)$ by setting 
\begin{align*}
    &\mathcal{N}_{L_t}(f) =0 \quad \text{for } f \in C^\infty(M), \\
    &(\mathcal{N}_{L_t} \alpha)(X,Y) = \alpha( (I+L_t)^{-1} N_{L_t} (X,Y)) \quad \text{for } \alpha \in A^1(M),
\end{align*}
and then extending it as a graded-symmetric derivation to the whole $A^\bullet(M)$.

To understand the action of the differential at time $t \neq 0$, we need the following result, cf.\ Lemma 3.4 in \cite{BT13}.

\begin{lemma}\label{lemma:commutator}
    Let $L_t \in \End (TM)$ and $\rho_t$ be given by \eqref{eq:rho}. Then
    \begin{equation}\label{eq:commutator}
    d \circ \rho_t = \rho_t \circ (d + [\tau_{L_t},d] - \mathcal{N}_{L_t}).
    \end{equation}
\end{lemma}
\begin{proof}
    The lemma is proved in the same way as Lemma 3.4 in \cite{BT13}. We just need to observe that the proof given in \cite{BT13} does not use the fact that the fundamental form $\omega$ is closed.
\end{proof}

We extend the definition of the operators $\partial$ and $\bar \partial$ to vector valued differential forms as follows. The operator $\bar \partial$ extends to the spaces $A^{0,q} \otimes TM$ as in Definition 2.10 in \cite{BM10}. Explicitly, if $X \in TM$, then $\bar \partial X \in A^{0,1}(M) \otimes TM$ is given by
\begin{equation}\label{eq:dbar:ex}
(\bar \partial X) (Y) \coloneqq \frac{1}{2} ( [Y, X] + J[JY, X]) - \frac{1}{4} N_J(X,Y),
\end{equation}
for all $Y \in TM$. If $\alpha \in A^{0,q} \otimes TM$, then $\bar \partial \alpha \in A^{0,q+1}(M) \otimes TM$ is given by
\begin{align*}
(\bar \partial \alpha) (X_0, \ldots, X_q) \coloneqq &\sum_{j=1}^q (-1)^j ( \bar \partial (\alpha(X_0, \ldots, \hat{X}_j, \ldots, X_q)))(X_j)\\
+ &\sum_{0 \le j \le k \le q} (-1)^{j+k} \alpha( [[X_j, X_k]], X_0, \ldots, \hat{X}_j, \ldots, \hat{X}_k, \ldots, X_q),
\end{align*}
where $[[X_j, X_k]]$ is given by
\[
[[X_j, X_k]] \coloneqq \frac{1}{2} ([X,Y] - [JX, JY]) + \frac{1}{4} N_J(X,Y).
\]

The operator $\partial$ can be extended to tensor fields $A^{0,q}(M) \otimes \Lambda^p (TM^{1,0})$ as follows. Fix a $J_0$ compatible metric on $M$, and consider the operator
\[
*_h \colon \Lambda^p (TM^{1,0}) \longrightarrow A^{m-p, 0}(M),
\]
defined by the relation
\[
\alpha \wedge *_h \beta = \langle \alpha | \beta \rangle \sigma_0,
\]
where $\alpha \in A^{p,0}(M)$, $\beta \in \Lambda^p (TM^{1,0})$, and $\langle \cdot | \cdot \rangle$ denotes the duality pairing induced by the metric. Then $\partial$ induces an operator
\[
\eth   \colon A^{0,q}(M) \otimes \Lambda^p (TM^{1,0}) \longrightarrow A^{0,q}(M) \otimes \Lambda^{p-1} (TM^{1,0}),
\]
given by $\eth \coloneqq (\id \otimes *_h)^{-1} \circ \partial \circ (\id \otimes *_h)$. The action of $\eth$ on $L \in \mathcal{L}$ can be described in terms of $\partial$ and $\tau_L$ as
\begin{equation}\label{eq:BTT}
\partial (\tau_L \sigma_0) = \eth   (L \wedge \sigma_0).
\end{equation}

Before proving the main result of this section, we state the following

\begin{lemma}\label{lemma:zero:square}
    Suppose that $J_0$ is integrable. Then the following identities hold:
    \[
    \bar \partial^2=0, \qquad \eth  ^2 =0, \qquad [\eth  , \bar \partial] =0.
    \]
\end{lemma}
\begin{proof}
    The only non-obvious identity is the first one. Its proof involves a long, although elementary, computation using the definitions and is omitted.
\end{proof}

We are ready to describe the tangent space to the moduli space of almost complex structures with PH trivial canonical bundle.

\begin{theorem}\label{thm:tangent:space}
    Let $M$ be a compact almost complex $2m$-manifold. Let $(J_0, \sigma_0) \in \mathcal{D}_{PH}$ and consider the moduli space
    \[
    \mathcal{M} \coloneqq  \mathcal{D}_{PH} / \Diff_0(M).
    \]
    If $J_0$ is integrable and it satisfies the $\partial \bar \partial$-lemma on $(m,1)$-forms, then we have that
    \[
    T_{(J_0, \sigma_0)} \mathcal{M} = \frac{\mathcal{L} \cap \ker \eth   \bar \partial}{\mathcal{L} \cap \Ima \bar \partial}
    \]  
\end{theorem}
\begin{proof}
    Fix a PH trivialization $\sigma_0$ of $K_{0}$, which is unique up to multiplication by non-zero constants, and let $(J_t, \sigma_t)$, $\abs{t} < \varepsilon$, be a smooth curve in $\mathcal{D}_{PH}$. By taking the projection on bidegree $(m,1)$ at time $t$ of equation \eqref{eq:commutator} applied to $F_t \sigma_0$, and by recalling that $\sigma_t = F_t \rho_t(\sigma_0)$, we have that
    \[
    \bar \partial_t \sigma_t = \rho_t \circ (\bar \partial + \pi^{m,1}_0 \circ [\tau_{L_t},d]) (F_t \sigma_0) + o(t),
    \]
    where $\bar \partial_t$ is the $\bar \partial$-operator at time $t$ and $\pi_0^{m,1}$ is the projection on bidegree $(m,1)$ at time $t=0$. One has that $\bar \partial_t \sigma_t =0$ if and only if
    \[
    \bar \partial  (F_t \sigma_0) + \pi^{m,1}_0 \circ [\tau_{L_t},d] (F_t \sigma_0) + o(t) =0.
    \]
    Taking time-derivative and evaluating at $t=0$, we get that
    \[
    \bar \partial(f \sigma_0) + \pi_0^{m,1} ( [\tau_L,d] \sigma_0 ) = 0,
    \]
    where $F_t = 1 + tf + o(t)$. By the initial condition $\bar \partial  \sigma_0 =0$ and the integrability of $J_0$, the equation reduces to
    \[
    \bar \partial  (f \sigma_0 ) = \partial  (\tau_L \sigma_0),
    \]
    and, by equation \eqref{eq:BTT}, we get the equality
    \[
    \bar \partial  (f \sigma_0 ) = \eth   (L \wedge \sigma_0).
    \]
    By applying $\bar \partial $ to both sides, we see that $\bar \partial  \eth   (L \wedge \sigma_0) = 0$. Conversely, if $L \in \mathcal{L}$ satisfies $\bar \partial  \eth   (L \wedge \sigma_0) = 0$, then $\bar \partial  \partial  (\tau_L \sigma_0) = 0$. Since the $\partial \bar \partial$-lemma holds for the $(m,1)$-form $\partial  (\tau_L \sigma_0)$, there exists $\alpha \in A^{m-1,0}(M)$ such that $\partial  \bar \partial  \alpha = \partial  (\tau_L \sigma_0)$, hence a smooth function $f$ such that $\bar \partial  (f \sigma_0 ) = \partial  (\tau_L \sigma_0)$. This proves that
    \[
    T_{(J_0,\sigma_0)} \mathcal{D}_{PH} \cong \{ (L,f) \in \mathcal{L} \times C^\infty (M): \bar \partial  (f \sigma_0) = \partial  (\tau_L \sigma_0) \} \cong \{ L \in \mathcal{L}: \eth   \bar \partial  L = 0 \}.
    \]

    To complete the description of the tangent space to $\mathcal{M}$, we have to quotient by the action of the group $\Diff_0 (M)$ of diffeomorphisms isotopic to the identity. Let $\phi_s$, for small $s \in \R$, be such a diffeomorphism, and let $X$ be the vector field generating it, that is, such that $\partial_{s=0} \phi^*_s = \mathcal{L}_X$. Under the action of $\phi_s$, the structure $J_0$ transforms as $d\phi_s^{-1} \circ J_0 \circ d\phi_s$. The corresponding infinitesimal action is given by
    \[
    \partial_{s=0}(d\phi_s^{-1} \circ J_0 \circ d\phi_s)(Y) = - \mathcal{L}_X (J_0Y) + J_0(\mathcal{L}_XY) = J_0[X,Y] -[X, J_0Y].
    \]
    Up to the action of $J_0$, we can identify the infinitesimal action with
    \[
    [X,Y] + J_0[X,JY] = (\bar \partial  X) (Y),
    \]
    where $\bar \partial  X$ is defined by \eqref{eq:dbar:ex}. Finally, since $\bar \partial ^2 =0$ by Lemma \ref{lemma:zero:square}, we have that $\eth   \bar \partial ^2 X = 0$, showing that $\Ima \bar \partial  \subseteq \ker \eth   \bar \partial $ and proving the inclusion $T_{(J_0, \sigma_0)} \mathcal{M} \subseteq (\mathcal{L} \cap \ker \eth   \bar \partial)/(\mathcal{L} \cap \Ima \bar \partial)$. The opposite inclusion follows by taking a first order deformation
    \[
    J_t = (I +tL) J_0(I+tL)^{-1},
    \]
    for every given $[L] \in (\mathcal{L} \cap \ker \eth   \bar \partial)/(\mathcal{L} \cap \Ima \bar \partial)$.
\end{proof}

\subsection{Tangent space to the space of structures with PH torsion canonical bundle.}

The results obtained in Section \ref{sec:PH:trivial} can be extended to the space of almost complex structures with PH torsion canonical bundle.
\vspace{.2cm}

Consider the space of almost complex structure on $M$ with PH torsion canonical bundle and torsion order bounded by $k$
\[
\mathcal{D}^k_{tor} = \{ (K, \eta) \in \J \times A^m(M): \eta \in A^{m,0}_K(M), \bar \partial_K (\eta^{\otimes l}) =0 \text{ for some $l \le k$}  \} \subset \mathcal{D}.
\]
Note that $\mathcal{D}^1_{tor} = \mathcal{D}_{PH}$. We have the following generalization of Theorem \ref{thm:tangent:space}.

\begin{theorem}\label{thm:PH:torsion}
    Let $M$ be a compact almost complex $2m$-manifold. Let $(J_0, \sigma_0) \in \mathcal{D}^k_{tor}$ and consider the moduli space
    \[
    \mathcal{M}^k_{tor} \coloneqq  \mathcal{D}^k_{tor} / \Diff_0(M).
    \]
    If $J_0$ is integrable and it satisfies the $\partial \bar \partial$-lemma on $(m,1)$-forms, then we have that
    \[
    T_{(J_0, \sigma_0)} \mathcal{M}^k_{tor} = \frac{\mathcal{L} \cap \ker \eth   \bar \partial }{\mathcal{L} \cap \Ima \bar \partial }
    \]  
\end{theorem}
\begin{proof}
    Fix $J_0 \in \mathcal{D}^k_{tor}$ satisfying the assumption of the theorem, and let $\Tilde{l}$ be the lowest common multiple of the torsion order of the canonical bundle in a neighborhood of $J_0$, which is well defined by the boundedness assumption. Let $J_t$, $\abs{t} < \varepsilon$, be a curve of almost complex structures in $\mathcal{D}^k_{tor}$ passing through $J_0$. Then the $\Tilde{l}$-th tensor power of the canonical bundle $K_t$ is PH trivial for all $t \in (-\varepsilon,\varepsilon)$, and determining the tangent space to $\mathcal{M}^k_{tor}$ is equivalent to determining the tangent space to the space of structures for which the canonical bundle is PH torsion with order that divides $\Tilde{l}$. 

    To do so, we slightly modify the proof of Theorem \ref{thm:tangent:space}. First, extend $\rho_t$ to a map between tensor powers of forms by setting $\rho_t (\sigma^{\otimes \Tilde{l}}) = (\rho_t \sigma)^{\otimes \Tilde{l}}$, where $\sigma$ is any section of $K_0$. Hence, any PH section $\sigma_t$ trivializing $K_t^{\otimes \Tilde{l}}$ can be written as
    \[
    \sigma_t = F_t (\rho_t \sigma_0)^{\otimes p},
    \]
    where $\Tilde{l} = p l_0$ and $l_0$ is the order of $K_0$. It is immediate to see that the claim of Lemma \ref{lemma:commutator} holds also for the operator $\rho_t$ extended to tensor powers, hence $\sigma_t$ must satisfy
    \[
    \bar \partial_t \sigma_t = \rho_t \circ (\bar \partial  + \pi^{m,1}_0 \circ [\tau_{L_t},d]) (F_t \sigma_0^{\otimes p}) + o(t).
    \]
    Thus, one has that $\bar \partial_t \sigma_t =0$ if and only if
    \[
    \bar \partial  (F_t \sigma_0^{\otimes p}) + \pi^{m,1}_0 \circ [\tau_{L_t},d] (F_t \sigma_0^{\otimes p}) + o(t) =0,
    \]
    that, at first order, gives the condition
    \[
    \bar \partial  (f \sigma_0^{\otimes p} ) - \partial  (\tau_L \sigma_0^{\otimes p})=0.
    \]
    We can develop the first term as
    \[
    \bar \partial  (f \sigma_0^{\otimes p} ) = \bar \partial f \otimes \sigma_0^{\otimes p} + p f (\bar \partial  \sigma_0) \otimes \sigma_0^{\otimes p-1} =  \bar \partial f \otimes \sigma_0^{\otimes p},
    \]
    where we used the fact that $\bar \partial  \sigma_0=0$, and the second one as
    \[
    \partial  (\tau_L \sigma_0^{\otimes p}) = p \partial  (\tau_L \sigma_0 \otimes \sigma_0^{\otimes p-1}) = p \partial  (\tau_L \sigma_0) \otimes \sigma_0^{\otimes p-1}.
    \]
    This implies that there exists $\hat{f} = \frac{f}{p}$ such that 
    \[
    \bar \partial  (\hat f \sigma_0 ) = \partial  (\tau_L \sigma_0).
    \]
    The rest of the proof proceeds as in Theorem \ref{thm:tangent:space} to get to the same conclusion.
\end{proof}

\subsection{Unobstructedness.}\label{section:unobstructedness}
We prove that, under suitable assumptions, the spaces $\mathcal{M}$ and $\mathcal{M}^k_{tor}$ are unobstructed.

\begin{theorem}\label{thm:unobstructed}
    Let $M$ be a compact almost complex $2m$-manifold. Let $J_0$ be a fixed almost complex structure on $M$ with PH torsion canonical bundle. Suppose that $J_0$ is integrable, that $H^{0,1}_{\bar \partial } =\{0\}$, and that the $\partial  \bar \partial $-lemma holds. Then the moduli space $\mathcal{M}^k_{tor}$ is unobstructed in a neighborhood of $J_0$.
\end{theorem}
\begin{proof}
    Recall that, by Theorem \ref{thm:PH:torsion}, the tangent space to $\mathcal{M}_{tor}^k$ at $(J_0, \sigma_0)$ can be identified with pairs 
    \[
    (L, f) \in \mathcal{L} \times C^\infty(M)^*
    \]
    satisfying the additional constraint $\bar \partial  (f \sigma_0) = \partial  (\tau_L \sigma_0)$. We first prove the theorem in the case $k=1$, corresponding to PH trivial canonical bundle.

    The proof of the theorem is a standard application of the implicit function theorem for the map $d \psi$, with $\psi$ given by
    \begin{align*}
    \psi \colon &\mathcal{L} \times C^\infty(M) \longrightarrow A^{0,1}(M)\\
    &(L, f) \longmapsto \iota_{\sigma_0} (\rho^{-1} \circ \bar \partial_K \circ \rho ( (1+f)\sigma_0)),
    \end{align*}
    where $K=(I+L)J_0(I+L)^{-1}$ and $\rho$ is given by \eqref{eq:rho} with $L_t=L$. By definition of $\psi$, the set $\psi^{-1}(0)$ can be identified with the space of almost complex structures with PH trivial canonical bundle. Furthermore, by the computation done in the proof of Theorem \ref{thm:tangent:space}, the differential of $\psi$ at $(0,0)$ is
    \[
    d \psi_{(0,0)} (L,f) = \bar \partial  f - \eth   L.
    \]
    
    Fix an arbitrary $J_0$-compatible metric and let $\Delta_{\bar \partial} \coloneqq \bar \partial \bar \partial^* + \bar \partial^* \bar \partial$ be the associated Dolbeault Laplacian. The assumption $H^{0,1}_{\bar \partial } = \{ 0 \}$ guarantees that $A^{0,1}(M) =  \Ima \Delta_{\bar \partial} \cap A^{0,1}(M)$, which gives the inclusion $\Ima d\psi_{(0,0)} \subseteq \Ima \Delta_{\bar \partial}$. We show that the opposite inclusion holds as well. Let $\alpha \in \Ima \Delta_{\bar \partial }$. By the standard Hodge decomposition for $\bar \partial$ on compact complex manifolds, we have that 
    \[
    \alpha = \bar \partial  f + \bar \partial ^* \theta,
    \]
    for some $f \in C^\infty (M)$ and $\theta \in A^{0,2}$. Thus, the $(0,2)$-form $\bar \partial  \alpha = \bar \partial  \bar \partial ^* \theta$ is $\bar \partial $-exact and $\eth  $-closed, by definition of $\eth  $. Validity of the $\partial \bar \partial$-lemma implies that the $\eth   \bar \partial $-lemma holds. Hence, there exists $L \in \mathcal{L}$ such that $\bar \partial  \bar \partial ^* \theta = \eth   \bar \partial  L$, giving the equality
    \[
    \alpha = \bar \partial  f - \eth   L
    \]
    and proving the opposite inclusion. Finally, we observe that the first-order operator $\psi(L,\cdot)$ is elliptic since, up to invertible maps, it behaves as the operator $\bar \partial_t$, which is elliptic, so that we can apply the implicit function theorem and standard regularity of elliptic operators to obtain that there exists a curve $(J_t,\sigma_t)$ passing through $(J_0,\sigma_0)$ whose tangent vector is $(L,f)$. This shows unobstructedness in the PH trivial case. For the PH torsion case, it is enough to work with a suitable tensor power of the canonical bundle as in the proof of Theorem \ref{thm:PH:torsion}.
\end{proof}

\subsection{Application to complex surfaces.}
We can combine the results of Theorems \ref{thm:iitaka}, \ref{thm:PH:torsion} and \ref{thm:unobstructed} to give a complete description of non-integrable small deformations of $K3$ and Enriques surfaces that have Kodaira dimension $0$. 

\begin{theorem}\label{thm:K3}
    Let $(M,J_0)$ be a $K3$ or Enriques surface and let $\mathcal{K}_0(M)$ be the space of almost complex structures on $M$ with vanishing Kodaira dimension. Then the moduli space $\mathcal{K}_0(M)/\Diff_0(M)$ is unobstructed in a neighborhood of $J_0$. Moreover, small deformations of $J_0$ with vanishing Kodaira dimension are in bijection with classes
    \[
    [L] \in  \frac{\mathcal{L} \cap \ker \eth   \bar \partial}{\mathcal{L} \cap \Ima \bar \partial}.
    \]
\end{theorem}
\begin{proof}
    We work locally in a neighborhood of $J_0 \in \J$. $K3$ and Enriques surfaces are K\"ahler surfaces with $b_1=0$, hence they satisfy all the assumptions of Theorems \ref{thm:PH:torsion} and \ref{thm:unobstructed}. This implies that $\mathcal{M}^k_{tor}$ is unobstructed in a neighborhood of $J_0$ for every fixed $k$, and that its tangent space is
    \[
    T_{(J_0, \sigma_0) } \mathcal{M}^k_{tor} = \frac{\mathcal{L} \cap \ker \eth   \bar \partial}{\mathcal{L} \cap \Ima \bar \partial},
    \]
    where $\sigma_0$ is a PH trivialization of a tensor power of the canonical bundle. For every $k$, we have the inclusion $\mathcal{M}^k_{tor} \subseteq \mathcal{K}_0(M)/\Diff_0(M)$. Conversely, let $J_1$ be a small deformation of $J_0$ such that $\kappa_{J_1}=0$. Then $J_1$ is still tamed and we have $c_1(J_1)= c_1(J_0)$, which is a torsion class, by the invariance of $c_1$ under small deformations. Hence, by Corollary \ref{cor:kodaira:4}, there exists $k$ such that $K_{J_1}$ is PH torsion of order $k$, giving the equality
    \[
    \mathcal{K}_0(M)/\Diff_0(M) = \bigcup_{k=1}^\infty \mathcal{M}^k_{tor}.
    \]
    The spaces $\mathcal{M}^k_{tor}$ are nested into each other and are all contained in $\J$, which is a tamed Frechet manifold. In addition, they have the same tangent space at $(J_0, \sigma_0)$ and are unobstructed. Hence, by the Nash--Moser implicit function theorem for tamed Frechet spaces, see, e.g., Theorem 51.17 in \cite{KM97}, they are locally diffeomorphic around $(J_0, \sigma_0)$. This shows that, in a neighborhood of $J_0$, we have that
    \[
    \mathcal{K}_0(M)/\Diff_0(M) = \mathcal{M}^{k_0}_{tor},
    \]
    where $k_0$ is the torsion order of $K_{J_0}$, and completes the proof of the theorem.
\end{proof}

\printbibliography

\end{document}